\documentclass[12pt]{amsart}
\usepackage{amsmath,amssymb,xypic,amsfonts,mathrsfs,subfigure,psfrag}
\usepackage{times,array,delarray,graphicx,dsfont,xspace,a4wide,epsfig}
\usepackage[numbers,sort&compress]{natbib}
\usepackage{hyperref}
\usepackage{rotating}
\graphicspath{{./}{./figs/}}

\theoremstyle{definition}
\newtheorem{definition}{Definition}[section]

\newtheorem{theorem}[definition]{Theorem}

\newtheorem{example}[definition]{Example}

\numberwithin{equation}{section}

\def\F{\mathbb{F}}

\def\Fy{\mathfrak{F}_Y}

\def\baralpha{\bar{\alpha}}
\def\barkappa{\bar{\kappa}}
\def\simkappa{\!\!\sim_{\kappa}}
\def\tsimkappa{\sim_{\kappa}}
\def\simbarkappa{\!\!\sim_{\bar\kappa}}
\def\tsimbarkappa{\sim_{\bar\kappa}}

\def\simalpha{\!\!\sim_{\alpha}}
\def\tsimalpha{\sim_{\alpha}}
\def\simbaralpha{\!\!\sim_{\bar\alpha}}
\def\tsimbaralpha{\sim_{\bar\alpha}}

\def\Acyc{\mathsf{Acyc}}

\def\<{\langle}
\def\>{\rangle}
\def\card#1{|#1|}
\def\Circle{\mathsf{Circ}}
\def\Circ{\mathsf{Circ}}
\def\vset{\mathrm{v}}

\def\varepsilon{\eta_e}

\def\Per{\mathsf{Per}}
\def\Aut{\mathsf{Aut}}

\def\nor{{\sf nor}\xspace }
\def\Nor{\mathsf{Nor}}

\let\to=\longrightarrow

\begin{document}
\bibliographystyle{amsplain}

\title{Update Sequence Stability in Graph Dynamical Systems}

\author{Matthew Macauley}
\address{Department of Mathematical Sciences \\
Clemson University \\
Clemson, SC 29634, USA}
\email{macaule@clemson.edu}

\author{Henning S.~Mortveit} \address{Department of Mathematics,
Virginia Tech, Blacksburg, VA 24061, USA} \email{henning@vt.edu}

\date{\today}

\begin{abstract}
  In this article, we study finite dynamical systems defined over
  graphs, where the functions are applied asynchronously. Our goal is
  to quantify and understand stability of the dynamics with respect to
  the update sequence, and to relate this to structural properties of
  the graph. We introduce and analyze three different notions of
  update sequence stability, each capturing different aspects of the
  dynamics. When compared to each other, these stability concepts
  yield vastly different conclusions regarding the relationship
  between stability and graph structure, painting a more complete
  picture of update sequence stability.
\end{abstract}


\maketitle

\section{Introduction}\label{sec:intro}

A graph dynamical system consists of $(i)$ an underlying graph $Y$
where each vertex has a state from a finite set $K$, $(ii)$ a sequence
of vertex functions that updates each vertex state as a function of
its neighbors, and $(iii)$ an update sequence that prescribes how to
compose the functions to get the global dynamical system map. Examples
of graph dynamical systems in the literature include systems with
synchronous update such as cellular automata~\cite{Wolfram:86a} and
Boolean networks~\cite{Shmulevich:02}, those with asynchronous
updates, e.g., sequential dynamical systems~\cite{Mortveit:07} and
functional linkage networks~\cite{Karaoz:04}, and many probabilistic
versions of the above systems. In this article, we focus on systems
with asynchronous schedules, with the broad goal of understanding the
sensitivity of the global dynamics with respect to changes in the
update sequence.

First, it is necessary to define what update order stability means,
and how to quantify such a concept. There are many possible answers to
this question, and there are vastly different ways to do this based on
existing literature. After a brief introduction of the main concepts
in Section~\ref{sec:prelims}, we devote the following three sections
to presentations of three different aspects of update sequence
stability, each focusing on different properties of the dynamics. In
Section~\ref{sec:equivalence}, we look at the structure of the phase
space and its robustness to changes in the update sequence. In
Section~\ref{sec:limitsets}, we analyze the number of attractor cycles
that can be reached from a given initial state under different
asynchronous update sequences. Finally, in
Section~\ref{sec:wordindependence}, we look at the collection of
periodic states and how they change (setwise) under different update
sequences.

Many of the concepts described above are scattered throughout the
literature, usually studied outside of the theme of update sequence
stability. This paper is the first of its kind to bring these ideas
together with this common goal in mind. It is unrealistic to expect to
develop a complete unifying theory of update order stability due to
the sheer size, diversity, and inherent randomness of complex
systems. However, there is still much be gained from such a study.
One of the themes of this paper is understanding how certain stability
properties are generally correlated with physical properties of the
underlying graph, such as edge density. For example, are systems over
sparse graphs or dense graphs more stable? Not surprisingly, different
notions of stability can lead to vastly different conclusions. We hope
that this paper will clarify many of these ideas and set the
groundwork for future investigations on the stability of asynchronous
finite dynamical systems.

\section{Preliminaries}\label{sec:prelims}

Let $Y$ be an undirected graph with vertex set
$\vset[Y]=\{1,\dots,n\}$. Each vertex $i$ has a state $x_i\in K$,
where $K$ is a finite set (frequently $K=\F_2=\{0,1\}$), and a
\emph{vertex function} $f_i$ that updates its state as a function of
the states of its neighbors (itself included). When applying these
functions asynchronously, it is convenient to encode $f_i$ as a
\emph{$Y$-local function} $F_i\colon K^n\to K^n$, where
\[
F_i(x_1,\dots,x_n)=(x_1,\dots,x_{i-1},f_i(x_1,\dots,x_n),x_{i+1},\dots,x_n)\;.
\]
In the remainder of this paper, we will focus on a class of graph
dynamical systems with asynchronous update schedules, called
\emph{sequential dynamical systems}~\cite{Mortveit:07}. These systems
have a prescribed update schedule dictated by a sequence of vertices
$\omega\in W_Y$, where $W_Y$ is the set of non-empty words over
$\vset[Y]$.
\begin{definition}
  A \emph{sequential dynamical system} (SDS) is a triple
  $(Y,\Fy,\omega)$, where $Y$ is an undirected graph, $\Fy$ a sequence
  of $Y$-local functions, and $\omega\in W_Y$ a length-$m$ update
  sequence. The \emph{SDS map} is the composition
  \[
  [\Fy,\omega]=F_{\omega_m}\circ F_{\omega_{m-1}}\circ\cdots\circ
  F_{\omega_2}\circ F_{\omega_1}\;.
  \]
\end{definition}
Henceforth, we will primarily restrict our attention to the case where
$\omega\in S_Y$, the set of permutations of $\vset[Y]$, i.e., words
where each vertex appears precisely once. To emphasize a permutation update
sequence, we will use $\pi$ instead of $\omega$. All of the concepts
carry over to general asynchronous systems under slight modification,
though not much additional insight is gained in doing so, whereas the
notation becomes more involved and the main concepts obscured.

The \emph{phase space} of a finite dynamical system $\phi\colon K^n\to
K^n$ is the directed graph $\Gamma(\phi)$ with vertex set $K^n$ and
edges $(x,\phi(x))$ for each $x\in K^n$. Clearly, $\Gamma(\phi)$
consists of disjoint cycles (periodic states), with directed trees
(transient states) attached.

\section{Equivalence of Dynamics}\label{sec:equivalence}

Any reasonable notion of stability of dynamics must be coupled with a
notion of equivalence of dynamics. For dynamical systems, a natural
place to begin is with the structure of the phase space, since this
encodes the global dynamics. We present here several different types of
equivalence, which vary in what is meant by having equivalent structure.
\begin{definition}
  Let $\phi,\psi\colon K^n\to K^n$ be finite dynamical systems. Then
  $\phi$ and $\psi$ are:
  \begin{enumerate}
  \item \emph{functionally equivalent} if $\phi=\psi$.
  \item \emph{dynamically equivalent} if $\phi\circ h=h\circ\psi$ for
    some bijection $h\colon K^n\to K^n$.
  \item \emph{cycle equivalent} if $\phi|_{\Per(\phi)}\circ
    h=h\circ\psi|_{\Per(\psi)}$ for some bijection
    $h\colon\Per(\phi)\to\Per(\psi)$.
  \end{enumerate}
\end{definition}
Clearly, $\phi$ and $\psi$ are functionally equivalent iff their phase
spaces are identical, dynamically equivalent iff their phase spaces
are isomorphic, and cycle equivalent iff their phase spaces are
isomorphic when restricted to the periodic states. For a fixed
sequence $\Fy$ of $Y$-local functions, different update sequences give
rise to different SDS maps, and each notion of equivalence mentioned
above partitions the set of all such maps into equivalence
classes. The number of SDS maps up to equivalence obtainable by
varying the update sequence is a measure of stability. As we will see,
each type of equivalence of dynamics corresponds with a combinatorial
equivalence relation on the set of update sequences $S_Y$ (or $W_Y$),
which in turn is strongly tied to the structure of the underlying
graph $Y$. Thus, knowledge of the structure of these equivalence
classes is central to quantifying update sequence stability, and its
relationship to the base graph. We will summarize this below.

\subsection{Functional Equivalence.}

Every update sequence $\pi\in S_Y$ defines a partial ordering $<_\pi$
of $\vset[Y]$, and this is naturally represented via the set
$\Acyc(Y)$ of acyclic orientations on $Y$. Explicitly, $\pi$ gives
rise to the acyclic orientation $O_\pi$ of $Y$ obtained by orienting
each edge $\{i,j\}$ of $Y$ as $(i,j)$ if $i$ occurs before $j$ in
$\pi$ and as $(j,i)$ otherwise. The partial order $<_\pi$ is defined
by $i<_\pi j$ iff there is a directed path from $i$ to $j$ in $O_Y$.
If $\pi,\sigma\in S_Y$ are linear extensions of the
same partial ordering, then $[\Fy,\pi]=[\Fy,\sigma]$, and we say that
$\pi\tsimalpha\sigma$. In fact, we  have a bijection
\[
\phi_\alpha \colon S_Y/\simalpha\longrightarrow\Acyc(Y)\;,
\]
hence $\alpha(Y):=|\Acyc(Y)|$ is an upper bound for the number of
distinct (permutation) SDS maps over $\Fy$ \emph{up to functional
equivalence}. This bound is sharp (see~\cite{Barrett:01a}), and it is
counted by the Tutte polynomial, as
$\alpha(Y)=T_Y(2,0)$. Additionally, the bijection $\phi_\alpha$ allows
us to construct a transversal of update sequences from $S_Y$, rather
than considering all $n!$ update sequences.

\subsection{Dynamical Equivalence.}

Let $\Aut(Y)$ denote the automorphism group of $Y$, and define the
$S_n$-action on $K^n$ by
$\sigma(x_1,\ldots,x_n)=(x_{\sigma^{-1}(1)},\ldots,x_{\sigma^{-1}(n)})$. The
group $\Aut(Y)$ acts on the set $S_Y/\simalpha$ by
$\gamma\cdot[\pi]_\alpha=[\gamma\pi]_\alpha$. If $\Fy$ is
$\Aut(Y)$-invariant (i.e., $\gamma\circ F_v=F_{\gamma(v)}\circ\gamma$
for each local function, and each $\gamma\in\Aut(Y)$), then two SDS
maps over update sequences in the same orbit are conjugate, by
$\gamma\circ [\Fy,\pi]\circ\gamma^{-1}=[\Fy,\gamma\pi]$
(see~\cite{Mortveit:07}). This puts an equivalence relation on $S_Y$
(and hence on $\Acyc(Y)$) that we denote by $\tsimbaralpha$. Since
conjugate maps are dynamically equivalent, the function
$\baralpha(Y):=|\Acyc(Y)/\simbaralpha\!|$ is an upper bound for the
number of SDS maps over $\Fy$ \emph{up to dynamical equivalence}. This
quantity, which counts the number of orbits in $\Acyc(Y)$ under
$\Aut(Y)$, can be computed though Burnside's lemma, as
\begin{equation*}
  \baralpha(Y)=\frac{1}{\card{\Aut(Y)}}\sum_{\gamma\in\Aut(Y)}
  \alpha(\<\gamma\>\setminus Y)\;.
\end{equation*}
Here, $\<\gamma\>\setminus Y$ denotes the \emph{orbit graph} of the
cyclic group $G=\<\gamma\>$ and $Y$,
(see~\cite{Barrett:01a,Barrett:03a}). We believe this bound is sharp,
and have shown this for special graph classes, but have no general
proof.

\subsection{Cycle Equivalence.}

To see how cycle equivalence arises in the SDS setting, we first
define a way to modify an acyclic orientation, called a
\emph{source-to-sink} operation, or a \emph{click}. This consists of
choosing a source vertex and reversing the orientations of all
incident edges. The reflexive transitive closure of this puts an
equivalence relation on $\Acyc(Y)$ (and on $S_Y$), denoted
$\tsimkappa$. Thus, two acyclic orientations related by a sequence of
clicks are $\kappa$-equivalent. In~\cite{Macauley:09a}, it is shown
that if $\pi\tsimkappa\sigma$, then $[\Fy,\pi]$ and $[\Fy,\sigma]$ are
cycle equivalent. Therefore, $\kappa(Y):=|\Acyc(Y)/\simkappa\!|$ is an
upper bound for the number of SDS maps over $\Fy$ \emph{up to cycle
equivalence}. This quantity is also counted by the Tutte polynomial,
by $\kappa(Y)=T_Y(1,0)$ (see~\cite{Macauley:08b}). As
with $\alpha$-equivalence, there is a nice transversal for
$\kappa$-equivalence, arising from the bijection
\begin{equation}
\label{eq:kappatrans}
  \phi_{\kappa,v} \colon \Acyc_v(Y) \longrightarrow \Acyc(Y)/\simkappa \;.
\end{equation}
where $\Acyc_v(Y)$ is the set of acyclic orientations of $Y$ where
vertex $v\in\vset[Y]$ is the unique source.

If additionally, $\Aut(Y)$ is nontrivial, we can use the action of
this group on $\Acyc(Y)/\simkappa$ to impose a coarser equivalence
relation $\tsimbarkappa$ on $\Acyc(Y)$. This is clear, since
dynamically equivalent SDS maps are trivially cycle equivalent. The
quantity $\barkappa(Y):=|\Acyc(Y)/\simbarkappa\!|$, which is the
number of orbits in $\Acyc(Y)/\simkappa$ under $\Aut(Y)$, is therefore
a stronger upper bound for the number of SDS maps over $\Fy$ up
to cycle equivalence when $\Fy$ is $\Aut(Y)$-invariant.

\subsection{Summary of Equivalence.}

To summarize, for each of the following notions of equivalence of
SDS maps, there is a corresponding equivalence relation
on $S_Y$ (and hence on $\Acyc(Y)$):
\begin{itemize}
  \item \emph{Functional equivalence}. Relation: $\tsimalpha$ with
    transversal $\Acyc(Y)$ via $\phi_\alpha$;
  \item \emph{Dynamical equivalence}. Relation: $\tsimbaralpha$ with
    transversal $\Acyc(Y)/\simbaralpha$;
  \item \emph{Cycle equivalence}. Relation: $\tsimkappa$ and
    transversal $\Acyc_v(Y)$ via $\phi_{\kappa,v}$ for $v\in\vset[Y]$.
\end{itemize}
In each of these three cases, update sequences that are equivalent
give rise to SDS maps with equivalent dynamics. Additionally, it is
clear that functionally equivalent systems are dynamically equivalent,
which in turn are cycle equivalent. In the general case where the
functions $\Fy$ are not $\Aut(Y)$-invariant, or when $\Aut(Y)=1$
(which is true of almost all random graphs~\cite{Bollobas:01}),
$\baralpha$-equivalence reduces down to $\alpha$-equivalence. The
difference between $\alpha$- and $\baralpha$-equivalence only captures
special symmetries in the graph. In light of this, we will focus on
$\alpha(Y)$ and $\kappa(Y)$ for our discussion of stability, though
$\baralpha(Y)$ was mentioned for sake of completeness.

Since the functions $\alpha(Y)$ and $\kappa(Y)$ bound the number of
SDS maps obtainable up to equivalence under different update
sequences, they can be considered as measures of stability. Computing
these quantities is in general, NP-hard, since they are evaluations of
the Tutte polynomial~\cite{Jaeger:90}. However, even though we can not
compute these functions exactly for large networks, an understanding
of the mathematics involved provides valuable insight about update
order stability. For example, consider two extreme cases: $(i)$ $Y$ is
an $n$-vertex tree, and $(ii)$ $Y=K_n$, the complete graph on $n$
vertices. In the first case, $\kappa(Y)=1$, and in the latter,
$\kappa(Y)=(n-1)!$.  Moreover, since $\alpha$ and $\kappa$ are
Tutte-Grothendieck invariants, if $Z\leq Y$ then
$\alpha(Z)\leq\alpha(Y)$ and $\kappa(Z)\leq\kappa(Y)$, i.e., $\alpha$
and $\kappa$ monotonically increase as additional edges are introduced
to the base graph. Two intermediate cases (between the tree and
complete graph) were studied in \cite{Macauley:09a} -- the graphs
corresponding to the radius-$1$ and radius-$2$ cellular automata
rules. For radius-$1$ rules (functions over the the circular graph
$\Circ_n$), there are at most $\kappa(Y)=O(n)$ cycle structures for
(permutation) SDS maps. In contrast, an SDS with a radius-$2$ rule has
base graph $Y=\Circ_{n,2}$, the circular graph where every vertex is
additionally connected to its distance-$2$ neighbors. In this case,
$\kappa(Y)=O(n\cdot 2^n)$. Though we cannot compute $\kappa(Y)$
explicitly for most graphs, we see a clear correlation between edge
density of the base graph and the stability of a dynamical system over
it. Thus, using $\alpha$ or $\kappa$-equivalence as a notion of
stability, we can make the general statement: \\
\[
\mbox{\emph{Update order stability and edge density are positively
correlated.}}
\]

\vskip5mm

\subsection{An Example.}

We conclude this section with an example illustrating the difference
between the various notions of equivalence described in this section.

\begin{example}
  \label{ex:running2}
  Let $Y=\Circ_4$, with vertex functions $\nor_3 \colon\F_2^3
  \longrightarrow\F_2$ where $\nor_3(0,0,0)=1$ and $\nor_3(x,y,z)=1$
  for all other arguments. Using update sequences $\pi = (1,2,3,4)$,
  $\pi'=(1,4,2,3)$ and $\pi'' = (1,3,2,4)$, the phase spaces
  are shown in Figure~\ref{fig:ex:norce}.
  \begin{figure}[ht]
    \centerline{
      \includegraphics[width=0.95\textwidth]{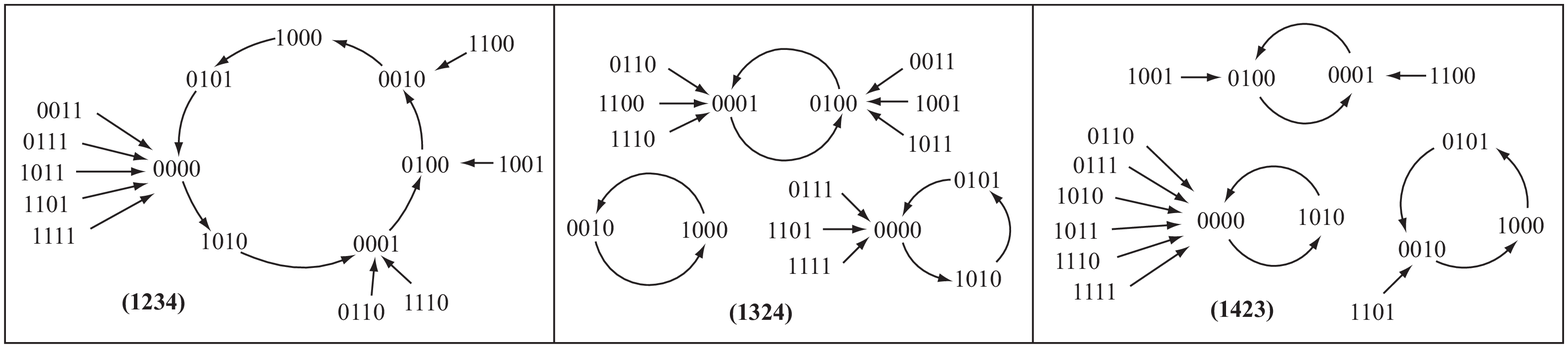}
    }
    \caption{The phase spaces of Example~\ref{ex:running2}.}
    \label{fig:ex:norce}
  \end{figure}
  None of these SDS maps are functionally or dynamically equivalent,
  but $[\Nor_Y,\pi']$ and $[\Nor_Y,\pi'']$ are cycle equivalent, as
  illustrated in the two rightmost phase spaces in
  Figure~\ref{fig:ex:norce}. In fact, it is elementary to compute
  $\bar\kappa(\Circle_4)=2$, thus the two cycle configurations in
  Figure~\ref{fig:ex:norce} are the \emph{only} possible such
  configurations up to isomorphism when $Y=\Circle_4$ and $K=\F_2$.
\end{example}

\section{Limit Sets and Reachability}\label{sec:limitsets}

Another natural measure for update sequence stability is a
quantification of how many attractor cycles can be reached from a
certain initial state. There are many ways to develop such a measure
-- one can compute the average number of periodic states reachable
from a random starting point, the maximum number of reachable periodic
states, the number of reachable attractor cycles, possibly normalized
for size, etc. For general complex networks, such computations are
intractable, and Monte Carlo type methods may be an
alternative. However, it is at times possible to estimate or bound
such quantities asymptotically. In this section, we will summarize
some recent results of one these measures that was developed to study
a the robustness of a gene annotation algorithm that basically uses a
threshold SDS~\cite{Karaoz:04}.

Threshold functions are common in discrete dynamical system
models~\cite{Kermack:27}. A \emph{Boolean threshold function} on $k$
arguments is a function $f_{k,m} \colon \F_2^k \longrightarrow \F_2$
such that $f(x_1,\ldots,x_k) = 1$ when at least $m$ of the arguments
$x_i$ are $1$, and $f(x_1,\ldots,x_k) = 0$ otherwise. An SDS is a
threshold SDS if every vertex function is a threshold function.
It is well-known that limit sets of threshold SDSs contain only fixed
points~\cite{Mortveit:07}, and additionally, that fixed points are
independent of update sequence.  This class plays a special role in
complexity theory since in many cases, threshold functions are the
most general function class for which problems are computationally
tractable~\cite{Barrett:03f,Barrett:01e}.

Functional linkage networks (FLNs) and the gene annotation algorithm
of~\cite{Karaoz:04} were among the motivations for considering the
type of update sequence stability measures discussed in this
section. An FLN (for a fixed gene ontology function $f$) is a graph
where the vertices represent proteins, and an edge is present if two
proteins are thought to share the same function, with an edge weight
$w_{ij}$ describing certainty. The vertices are assigned states from
$K=\{1,-1,0\}$ denoting whether they are annotated with $f$ (state
$1$), are not annotated with $f$ (state $-1$), or if annotation is
unknown (state $0$). Each vertex $i$ is updated by changing its state
to $1$ if $\sum w_{ij} x_j\geq T$, where the sum is taken over all
neighboring vertices, and $T$ is some fixed threshold, and to $-1$
otherwise. The algorithm updates the vertices in some fixed
order chosen at random, and terminates when a fixed point is reached.

The FLN algorithm uses an asynchronous update schedule to ensure that
a fixed point is reached (there may be periodic cycles of length $>1$
under a synchronous update schedule). In light of this convenient
choice of update scheme, and the random choice of update sequence, one
would hope that such a gene annotation algorithm would exhibit update
order stability. In other words, the final (fixed point) state reached
should be robust with respect to changes in update sequence and
perturbations of the initial state, in that the fixed point reach has
little or no dependence on the update sequence. Unfortunately, such
stability cannot always be guaranteed. As an example, there are
general classes of threshold SDSs sharing many of the essential FLN
algorithm properties, but which exhibit exponential update sequence
instability~\cite{Kumar:09} for many choices of the initial state. Not
surprisingly, such instability strongly depends on the choice of
underlying graph.

The {\em $\omega$-limit set} of $\phi\colon K^n\to K^n$ from $x\in
K^n$ is the set $\omega(\phi;x)$ of periodic states $x\in K^n$ such
that $\phi^m(x)=y$ for some $m\geq 0$. This is a specialization of the
classical definition of $\omega$-limit set to the case of finite phase
spaces. If $\mathcal{P}\subset S_Y$ is a collection of update
sequences, then the $\omega$-limit set of $\Fy$ from $x$ with respect
to $\mathcal{P}$ is defined as
\begin{equation*}
  \omega_{\mathcal{P}}(x)=\bigcup_{\pi\in\mathcal{P}}\omega([\Fy,\pi];x)\;.
\end{equation*}
If there are large periodic cycles in the phase space of the SDS, then
the size of an $\omega$-limit set may not necessarily be too
insightful. However, if there only fixed points, as is common in many
applications, then $\omega_{\mathcal{P}}(x)$ is a direct measure of
stability that describes how many fixed points can be reached from $x$
by choosing the update schedule from $\mathcal{P}$. For a set of
functions $\Fy$, define
\[
\omega(\Fy)=\max\left|\{\omega_{S_Y}(x)\mid x\in K^n\}\right|\;.
\]
This counts the maximum possible number of periodic points that can be
reached from any state by variation of the update order. As in the
previous section, we begin by considering two extreme cases; when the
base graph is a tree, and when it is the complete graph $K_n$. The
following result is proven in~\cite{Kumar:09}.
\begin{theorem}
  Let $\Fy$ be a sequence of 2-threshold functions. If $Y$ is the
  star-graph on $n$ vertices, then $\omega(\Fy)=2^n-n$. If $Y=K_n$,
  then $\omega(\Fy)=n+1$.
\end{theorem}
Additionally, this is extended to the random graph model $G(n,p)$ for
ranges of $p$ relevant to many applications.  The main result is that
for sparse graphs, $\omega(\Fy)=\Theta(2^n)$ (with high probability),
but for dense graphs, $\omega(\Fy)=\Theta(n)$. While these systems are
deliberately constructed to possess dynamical instability, they
nonetheless offer valuable insight into possible dynamics of such
systems. With this notion of stability with respect to update
sequence, one can make the following general statement, which is in
direct contrast to the conclusion of the previous section: \\
\[
\mbox{\emph{Update order stability and edge density are negatively
  correlated.}}
\]

\vskip5mm

\section{Periodic Points and Word Independence}\label{sec:wordindependence}

We come at last to our third and final stability concept for
asynchronous dynamical systems, which stems from recent work on word
independent systems~\cite{Hansson:05b, Macauley:08a,
Macauley:08e}. Like cycle equivalence, this considers invariance
properties of the periodic points of SDSs, rather than the transient
states. A sequence of $Y$-local functions $\Fy$ is \emph{word
independent} if the set of periodic states $\Per([\Fy,\omega])$ is the
same for all \emph{complete words} $\omega\in W_Y$ (words over
$\vset[Y]$ where each vertex occurs at least once). Note that word
independence only considers the periodic states as a set, and ignores
the orbit structure.  Word independence was first studied
in~\cite{Hansson:05b}, where it was shown that most classic symmetric
functions such as the logical AND, OR, NAND, NOR functions, along with
their sums, were word independent, for asynchronous cellular automata
(i.e., over the circular graph). Other common functions, such as
parity, minority, and majority functions had this property as
well. Recently, in a series of two papers~\cite{Macauley:08a,
Macauley:08e}, we have expanded this study to all $256$ elementary
cellular automata rules, and classified the word independent ones,
which are independent of the size of the underlying (circular) graph.
\begin{theorem}
  Precisely $104$ of the $256$ elementary cellular automata rules
  induce a word independent SDS over $\Circ_n$, for all $n>3$.
\end{theorem}
The classification of these $104$ rules is contained
in~\cite{Macauley:08a}, and the structure is further developed
in~\cite{Macauley:08e} with an analysis of their \emph{dynamics
groups}, which describe how the periodic states are permuted within
these systems as the local functions are applied. The latter paper
additionally contains more general results on dynamics groups. For
example, when the state space is $K=\{0,1\}$, as is typical in finite
dynamical systems research, the dynamics group is the homomorphic
image of a Coxeter group~\cite{Bjorner:05}, and for general $K$, an
Artin group.

The papers above are only the beginning of this line of work. There is
much to be done to better understand systems that fail to be word
independent. Invariance of periodic states under different update
sequences can be used as a notion of update sequence stability. This
can be further quantified as follows. When a sequence of functions is
not word independent, we can define a measure that captures how
``close'' it is to being word independent, by defining
\[
\rho(\Fy)=\frac{|\,x\in\Per[\Fy,\pi]\mbox{ for all }\pi\in S_Y|}
{|\,x\in\Per[\Fy,\pi]\mbox{ for some }\pi\in S_Y|}\;.
\]
Clearly, $0<\rho(\Fy)\leq 1$, with $\rho(\Fy)=1$ iff $\Fy$ is word
independent. This quantity measures how stable the periodic point sets
are as a whole, to changes in the update sequence. Under this
definition, word independent sequences of functions are the most
stable. In general, it is difficult to determine if a sequence of
functions is word independent. Though there may exist a correlation
between $\rho(\Fy)$ and the topology (e.g., edge density) of $Y$,
there is no good reason to believe that such a link should
exist. Rather, there are many examples of functions that are word
independent for all graphs, such as threshold functions, monotone
functions, functions with only fixed points, parity functions, and the
logical NOR functions. Therefore, it is not far-fetched to make the
following hypothesis: \\
\[
\mbox{\emph{There is little, if any, correlation between update order
stability and edge density.}}
\]
\vskip5mm

\section{Concluding Remarks}\label{sec:conclusion}

This article should raise a red flag for anyone anyone attempting to
study and characterize stability with respect to update sequence
variation. We saw how different notions of stability can lead to
different conclusions of the relationship to graph structure. In
summary, we looked at three properties of the dynamics and how robust
they are to changes in the update sequence. In particular, we
considered:
\begin{itemize}
\item The variation of structural properties of the phase space as the
  update sequence changes;
\item The variation of reachable limit sets as the update sequence changes;
\item The variation of the sets of periodic states as the update
  sequence changes.
\end{itemize}
For each of these notions, we defined natural measures quantifying
stability, and then inquired about the relationship between these
measures and the structure of the underlying graph. Though such a
question is quite open-ended, we began by looking at the extreme cases
of sparse graphs (e.g., trees) and dense graphs (e.g., complete
graphs) to gain some insight, before considering intermediate
cases. Up until now, these three notions of stability have only been
studied independently, and moreover, outside of the setting of
stability. In this article, we see first-hand how differently they
behave with respect to a basic property such as edge density of the
base graph. Perhaps there are ways to draw connections between them to
paint a more complete picture of stability.  It is our hope that this
article will help lay the groundwork for such future research.

\providecommand{\bysame}{\leavevmode\hbox to3em{\hrulefill}\thinspace}
\providecommand{\MR}{\relax\ifhmode\unskip\space\fi MR }
\providecommand{\MRhref}[2]{%
  \href{http://www.ams.org/mathscinet-getitem?mr=#1}{#2}
}
\providecommand{\href}[2]{#2}

\end{document}